\documentclass[12pt,oneside,draft]{amsart}
\usepackage{amssymb,amsmath}
\usepackage{graphpap}

 
\newtheorem{theorem}[subsection]{Theorem}
\newtheorem*{theorem*}{Theorem}

\newtheorem{proposition}[subsection]{Proposition}
\newtheorem{corollary}[subsection]{Corollary}

\theoremstyle{definition}
\newtheorem{definition}[subsection]{Definition}

\newtheorem{example}[subsection]{Example}
\theoremstyle{remark}
\newtheorem{remark}[subsection]{Remark}

\newcommand{\mt}[1]{\operatorname{#1}}
\newcommand{\EEE}{{\mathbb E}}
\newcommand{\DDD}{{\mathbb D}}
\newcommand{\AAA}{{\mathbb A}}
\newcommand{\QQ}{{\mathbb Q}}
\newcommand{\ZZ}{{\mathbb Z}}
\newcommand{\CC}{{\mathbb C}}
\newcommand{\OO}{{\mathcal O}}

\newcommand{\PP}{{\mathbb P}}

\newcommand{\TT}{{\mathcal T}}

\newcommand{\Sing}{\mt{Sing}}

\newcommand{\rank}{\mt{rank}}

\newcommand{\Exc}{\mt{Exc}}

\newcommand{\Spec}{\mt{Spec}}
\newcommand{\Gal}{\mt{Gal}}

\newcommand{\down}[1]{\llcorner #1 \lrcorner}

\newenvironment{outline}{\begin{proof}[The main steps of proof]}{\end{proof}}

\title{Classification of log Enriques surfaces with {\large $\delta=2$}}
\author{S.~A.~Kudryavtsev}

\date{}

\address{Department of Algebra, Faculty of Mathematics,
Moscow State Lomonosov University, 117234 Moscow, Russia}

\email{kudryav@mech.math.msu.su}

\begin{document}
\begin{abstract} Log Enriques surfaces with $\delta=2$ are classified.
\end{abstract}
\maketitle

\section*{\bf {Introduction}}
Let $S$ be a projective surface with klt singularities and with numerically trivial
canonical class $K_S\equiv 0$. Then $S$ is called {\it a log Enriques surface}.
\par
Let us consider the following invariant
\begin{eqnarray*}
\delta(S)=\#\Big\{E| E \ \text{is exceptional divisor with discrepancy}\\
a(E,0)\le -\frac67   \Big\}.
\end{eqnarray*}

By theorem \cite[5.1]{Sh2} we have $0\le \delta(S) \le 2$.
In this paper the classification of such surfaces with $\delta=2$ is given
(see theorems \ref{main1} and \ref{main2}).
\par
The log Enriques surfaces often appear in many problems.
Such of them are the following ones:
the study of surface degeneration, the study of $K3$ surfaces
(in particular, see the last section of \S 2), the study of
Calabi--Yao varieties and the problem of inductive classification of strictly log
canonical singularities.
\par
In the latter problem the log Enriques surfaces can be realized as the exceptional
divisors of purely log terminal blow-ups of three-dimensional strictly
log canonical singularities.
For example, in the case of hypersurface singularities we have $K_S\sim 0$ and $S$ is a
$K3$ surface with Du Val singularities or $S$ is a surface obtained by
the contraction of section of elliptic surface \cite{Kud3}.
\par
The method applied in this paper was developed by
V.V.~Shokurov in
\cite{Sh2}. It works in any dimension for the classification of the varieties with
Kodaira dimension
$-\infty$ or $0$, extremal contractions and singularities.
\par
The number of different log Enriques surfaces with
$\delta=2$ is more then 1000 and hence (for example) to enumerate
their minimal resolutions is difficult. This approach allows to obtain
the short and complete classification.
\par
The other methods of log Enriques surface study were given in
\cite{Bl}, \cite{Z} and \cite{Z2}. In many papers the automorphisms of
$K3$ surfaces were investigated and hence some results about the structure
of log Enriques surfaces were obtained (for example, see \cite{OZ1}, \cite{OZ2}).

\par
In this paper we can also see that the classification of the "models"\ of
exceptional log Del Pezzo surfaces
(see definitions of exceptionality and model in \cite{Sh2}, see also
theorem \ref{ldmain}) implies the classification of log Enriques surface.
\par
The research was partially supported by a grant 02-01-00441 from
the Russian Foundation of Basic Research and a grant INTAS-OPEN
2000\#269.

\section{\bf {Preliminary facts and results}}

All varieties are algebraic and are assumed to be defined over
$\CC$, the complex number field.
The main definitions, terminology and notations used in the paper are
given in \cite{Koetal}, \cite{PrLect}.

\begin{definition}
Let $X$ be a normal variety and let $D=S+B$ be a subboundary on
$X$ such that $B$ and $S$ have no common components, $S$ is an
effective integral divisor and $\down{B}\le 0$. Then we say that
$K_X+D$ is \textit{$n$-complementary} if there is a $\QQ$-divisor
$D^+$ such that
\begin{enumerate}
\item
$n(K_X+D^+)\sim 0$ (in particular, $nD^+$ is an integral divisor);
\item
$K_X+D^+$ is lc;
\item
$nD^+\ge nS+\down{(n+1)B}$.
\end{enumerate}
In this situation the
\textit{$n$-complement} of
$K_X+D$ is $K_X+D^+$. The divisor $D^+$ is called
\textit{$n$-complement} too.
\end{definition}

\begin{theorem}\label{contr} \cite[3.1]{Sh2}
Let $(X/Z\ni P,D)$ be a log surface of local type
(i.e. $Z$ is not point), where $f\colon X\to Z\ni P$ is a contraction.
Assume that $-(K_X+D)$ is $f$-nef and $K_X+D$ is lc.
Then there exists an 1,2,3,4 or 6-complement of $K_X+D$ near
$f^{-1}(P)$.
\end{theorem}

\begin{theorem}\label{ldmain} \cite[\S 5]{Sh2}
Let $S$ be a log Enriques surface with $\delta=2$. Denote the exceptional
curves with discrepancy $a(\widetilde{C}_i,0)\le -\frac67$ by
$\widetilde{C}_1$ and
$\widetilde{C}_2$. Let $f\colon \widetilde S\to  S$ be an extraction of
$\widetilde{C}_i$
$($i.e. $f$ is a birational contraction such that
$\Exc f=\widetilde C_1\cup \widetilde C_2)$.
\\
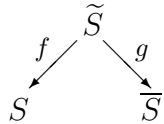
\begin{center}
\begin{picture}(80,46)(0,0)
\put(14,10){$S$}
\put(40,40){\vector(-1,-1){18}}
\put(42,42){$\widetilde S$}
\put(50,40){\vector(1,-1){18}}
\put(64,10){$\overline{S}$}
\put(24,34){\footnotesize{$f$}}
\put(62,34){\footnotesize{$g$}}
\end{picture}
\end{center}
\par
Then there exists a birational contraction
$g\colon \widetilde S\to  \overline{S}$ with the following properties:
$\rho(\overline{S})=1$, $g$
doesn't contract the curves $\widetilde{C}_i$. Put $C_i=g(\widetilde{C}_i)$.
Moreover, the pair
$\big(\overline{S},a(\widetilde{C}_1,0)C_1+a(\widetilde{C}_2,0)C_2\big)$ is one of the
following ones:
\begin{enumerate}
\item[$\mathrm{(A_2^6)}$]
$\overline{S}=\PP(1,2,3)$, $a(\widetilde{C}_1,0)=a(\widetilde{C}_2,0)=-\frac67$,
$C_1$ is a line
$\{x_1=0\}$, $C_2\in |-K_X|$. The curve $C_2$ must have an ordinary double point.
\item[$\mathrm{(I_2^2)}$]
$\overline{S}=\PP(1,2,3)$, $a(\widetilde{C}_1,0)=a(\widetilde{C}_2,0)=-\frac67$,
$C_1=\{x_3=0\}$,
$C_2=\{x_2^2+a_1x_1^4+a_2x_1^2x_2+x_1x_3=0\}$, $a_1,
a_2\in\CC$, $(a_1,a_2)\ne (0,0)$.
\end{enumerate}
\begin{outline} The divisor $K_{\widetilde{S}}-a(\widetilde{C}_1,0) \widetilde{C}_1-
a(\widetilde{C}_2,0) \widetilde{C}_2$ is klt. Therefore by theorem
\ref{contr} there exists an 1,2,3,4 or 6-complement near any point of
$\widetilde{S}$. Hence
$K_{\widetilde{S}}+\widetilde{C}_1+\widetilde{C}_2$ is lc by the definition
of complement. It is clear that
$-(K_{\widetilde{S}}+\widetilde{C}_1+\widetilde{C}_2)$ is not nef.
\par
Then there exists an extremal ray $R$ such that
$(K_{\widetilde{S}}+\widetilde{C}_1+\widetilde{C}_2)\cdot R>0$.
The contraction of this ray is birational and the curves
$\widetilde{C}_i$ are not contracted.
To prove these statements the theorem \ref{contr} is used.
We get a required model
$\overline{S}$ by repeating such procedure.
\par
Then such models are classified.
Taking into account the condition
$-K_{\overline{S}}+a(\widetilde{C}_1,0)C_1+
a(\widetilde{C}_2,0)C_2\equiv 0$
we obtain only two cases for $\overline{S}$.
\par
Since klt singularities are rational then the curves
$C_i$ must be rational. Hence in the case
$\mathrm{(A_2^6)}$ the curve $C_2$ must have an ordinary double point.
\end{outline}
\end{theorem}

\begin{definition} Let $S$ be a log Enriques surface. The minimal index of
complementary $I$ is called {\it a canonical index} of $S$, i.e.
$I=\min\{n\in\ZZ_{>0}\mid nK_S\sim 0\}$.
It is known that $I\le 21$ \cite{Bl}, \cite{Z}.
\end{definition}

\begin{corollary} \label{index}
Let $S$ be a log Enriques surface with $\delta=2$. Then
$I=7$.
\end{corollary}

So, the problem of classification of log Enriques surface with $\delta=2$
is to describe the following procedures.
At first we consider the extraction
$\widetilde S \to \overline{S}$ such that every exceptional divisor
$E$ has the discrepancy $a(E,\frac67C_1+\frac67C_2)=0$.
Then we contract the proper transforms of
$C_1$ and $C_2$.
The number of such procedures is finite by the following easy fact.
\par
Let $(X,D)$ be a klt pair. Then the number of divisors
$E$ of the function field $\mathcal{K}(X)$ with
$a(E,D)\le 0$ is finite \cite[lemma 3.1.9]{PrLect}.

\section{\bf {Classification of log Enriques surfaces with $\delta=2$ }}

\begin{proposition}\label{extr1}
\       \newline
\begin{enumerate}
\item Let us consider the pair $(X,\frac67C)\simeq(\CC^2, \frac67\{x=0
\})/\ZZ_2(1, 1)$. Then the extraction of all exceptional curves with discrepancies
$0$ is shown in the following figure.
\\
\normalfont
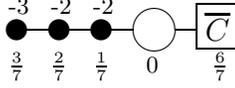
\begin{center}
\begin{picture}(100,30)(0,0)
\put(5,20){\circle*{8}}
\put(2,26){\scriptsize{-3}}
\put(2,4){\scriptsize{$\frac37$}}
\put(9,20){\line(1,0){8}}
\put(21,20){\circle*{8}}
\put(18,26){\scriptsize{-2}}
\put(18,4){\scriptsize{$\frac27$}}
\put(25,20){\line(1,0){8}}
\put(37,20){\circle*{8}}
\put(34,26){\scriptsize{-2}}
\put(34,4){\scriptsize{$\frac17$}}
\put(41,20){\line(1,0){8}}
\put(57,20){\circle{16}}
\put(65,20){\line(1,0){8}}
\put(54,4){\scriptsize{0}}
\put(73,16){\fbox{$\overline{C}$}}
\put(79,4){\scriptsize{$\frac67$}}
\end{picture}
\end{center}

{\it The proper transform of $C$ is denoted by $\overline{C}$. The numbers over
vertexes denote the self-intersection indexes of corresponding curves.
The numbers under vertexes
denote the discrepancies with opposite sign of corresponding curves.
The empty circle denotes the required exceptional curve with discrepancy $0$.
Its self-intersection index is always equal to $-1$ (on a minimal resolution).
Note that
$\overline{C}^2=C^2-\frac72$.
\item Let us consider the pair $(X,\frac67C)\simeq(\CC^2, \frac67\{x=0\})/\ZZ_3(1, 2)$.
Then the extraction of all exceptional curves with discrepancies
$0$ is shown in the following figure.}
\\
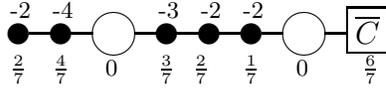
\begin{center}
\begin{picture}(150,30)(0,0)
\put(5,20){\circle*{8}}
\put(2,26){\scriptsize{-2}}
\put(2,4){\scriptsize{$\frac27$}}
\put(9,20){\line(1,0){8}}
\put(21,20){\circle*{8}}
\put(18,26){\scriptsize{-4}}
\put(18,4){\scriptsize{$\frac47$}}
\put(25,20){\line(1,0){8}}
\put(41,20){\circle{16}}
\put(38,4){\scriptsize{0}}
\put(49,20){\line(1,0){8}}
\put(61,20){\circle*{8}}
\put(58,26){\scriptsize{-3}}
\put(58,4){\scriptsize{$\frac37$}}
\put(65,20){\line(1,0){8}}
\put(77,20){\circle*{8}}
\put(74,26){\scriptsize{-2}}
\put(71,4){\scriptsize{$\frac27$}}
\put(81,20){\line(1,0){8}}
\put(93,20){\circle*{8}}
\put(90,26){\scriptsize{-2}}
\put(90,4){\scriptsize{$\frac17$}}
\put(97,20){\line(1,0){8}}
\put(113,20){\circle{16}}
\put(121,20){\line(1,0){8}}
\put(110,4){\scriptsize{0}}
\put(129,16){\fbox{$\overline{C}$}}
\put(135,4){\scriptsize{$\frac67$}}
\end{picture}
\end{center}
{\it The notations are as in the point $(1)$. Note that
$\overline{C}^2=C^2-\frac{14}3$.

\item Consider the pair $(X,\frac67C_1+\frac67C_2)\simeq(\CC^2, \frac67\{xy=0\})$.
Then the extraction of all exceptional curves with discrepancies
$0$ is shown in the following figure.}
\\
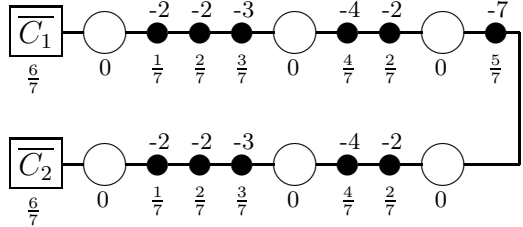
\begin{center}
\begin{picture}(200,90)(0,0)
\put(5,76){\fbox{$\overline{C_1}$}}
\put(11,60){\scriptsize{$\frac67$}}
\put(25,80){\line(1,0){8}}
\put(41,80){\circle{16}}
\put(39,64){\scriptsize{0}}
\put(49,80){\line(1,0){8}}
\put(61,80){\circle*{8}}
\put(58,64){\scriptsize{$\frac17$}}
\put(58,86){\scriptsize{-2}}
\put(65,80){\line(1,0){8}}
\put(77,80){\circle*{8}}
\put(74,64){\scriptsize{$\frac27$}}
\put(74,86){\scriptsize{-2}}
\put(81,80){\line(1,0){8}}
\put(93,80){\circle*{8}}
\put(90,64){\scriptsize{$\frac37$}}
\put(90,86){\scriptsize{-3}}
\put(97,80){\line(1,0){8}}
\put(113,80){\circle{16}}
\put(110,64){\scriptsize{0}}
\put(121,80){\line(1,0){8}}
\put(133,80){\circle*{8}}
\put(130,64){\scriptsize{$\frac47$}}
\put(130,86){\scriptsize{-4}}
\put(137,80){\line(1,0){8}}
\put(149,80){\circle*{8}}
\put(146,64){\scriptsize{$\frac27$}}
\put(146,86){\scriptsize{-2}}
\put(153,80){\line(1,0){8}}
\put(169,80){\circle{16}}
\put(166,64){\scriptsize{0}}
\put(177,80){\line(1,0){8}}
\put(189,80){\circle*{8}}
\put(186,64){\scriptsize{$\frac57$}}
\put(186,86){\scriptsize{-7}}
\put(193,80){\line(1,0){5}}
\put(198,80){\line(0,-1){50}}
\put(5,26){\fbox{$\overline{C_2}$}}
\put(11,10){\scriptsize{$\frac67$}}
\put(25,30){\line(1,0){8}}
\put(41,30){\circle{16}}
\put(38,14){\scriptsize{0}}
\put(49,30){\line(1,0){8}}
\put(61,30){\circle*{8}}
\put(58,14){\scriptsize{$\frac17$}}
\put(58,36){\scriptsize{-2}}
\put(65,30){\line(1,0){8}}
\put(77,30){\circle*{8}}
\put(74,14){\scriptsize{$\frac27$}}
\put(74,36){\scriptsize{-2}}
\put(81,30){\line(1,0){8}}
\put(93,30){\circle*{8}}
\put(90,14){\scriptsize{$\frac37$}}
\put(90,36){\scriptsize{-3}}
\put(97,30){\line(1,0){8}}
\put(113,30){\circle{16}}
\put(110,14){\scriptsize{0}}
\put(121,30){\line(1,0){8}}
\put(133,30){\circle*{8}}
\put(130,14){\scriptsize{$\frac47$}}
\put(130,36){\scriptsize{-4}}
\put(137,30){\line(1,0){8}}
\put(149,30){\circle*{8}}
\put(146,14){\scriptsize{$\frac27$}}
\put(146,36){\scriptsize{-2}}
\put(153,30){\line(1,0){8}}
\put(169,30){\circle{16}}
\put(166,14){\scriptsize{0}}
\put(177,30){\line(1,0){21}}
\end{picture}
\end{center}
{\it The notations are as in the point $(1)$. Notice that
$\overline{C}_i^2=C_i^2-6$ if $C_1+C_2$ is a (globally) reducible curve and
$\overline{C}^2=C^2-14$ if $C=C_1+C_2$ is an (globally) irreducible curve
with an ordinary double point.}
\end{enumerate}
\begin{proof} In the cases (1) and (2) at first
we consider the minimal resolution of singularity
$f\colon Y\to X$. Then $K_Y+\sum d_iD_i=f^*(K_X+\frac67C)$, where $D_i$ are the
irreducible divisors.
If $D_i\cdot D_j\ne 0$ for some $i\ne j$ and $d_i+d_j\ge 1$ then
let us blow-up their intersection point.
After the finite number of such blow-ups we will extract  all exceptional
divisors with discrepancies 0.
Note also that the exceptional divisor with discrepancy 0 appears if and only if
$d_i+d_j=1$.
\end{proof}
\end{proposition}

This proposition trivially implies the next two corollaries.
\begin{corollary}\label{extr2}
Let us consider the case $\mathrm{(A_2^6)}$ of theorem \ref{ldmain}.
Then the extraction of all exceptional curves with discrepancies
$0$ for the pair $(\overline{S},\frac67C_1+\frac67C_2)$
is shown in the following figure.
\end{corollary}
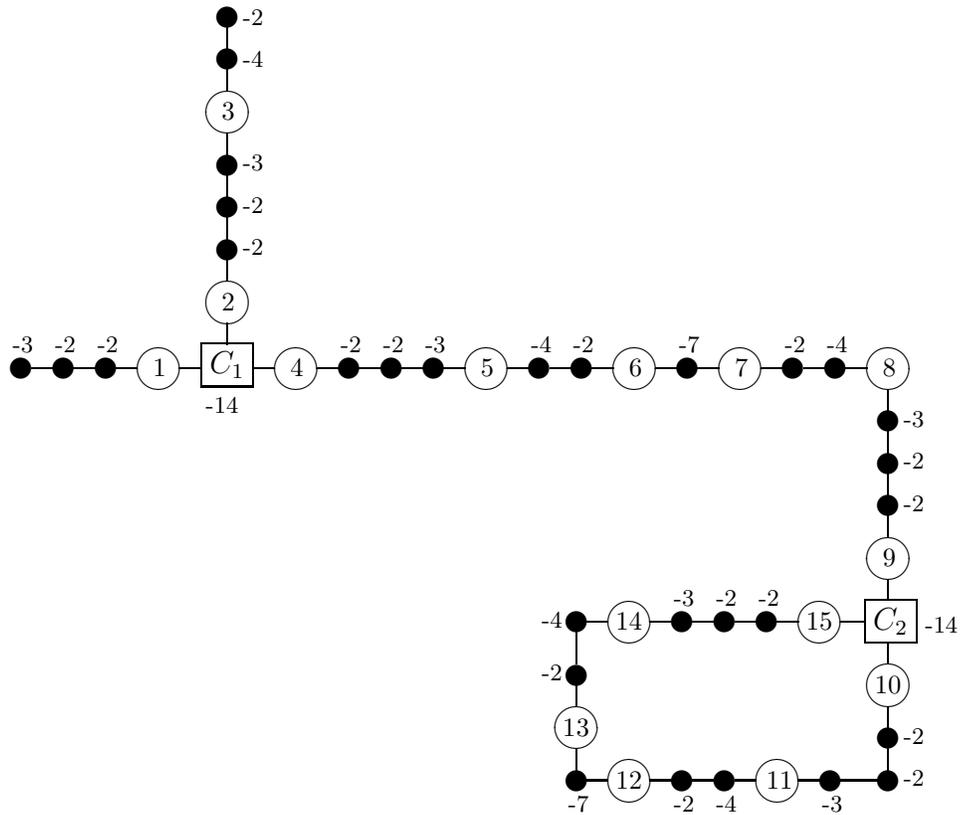
\begin{center}
\begin{equation*}
\begin{picture}(344,299)(0,-25)
\put(4,150){\circle*{8}}
\put(1,156){\scriptsize{-3}}
\put(8,150){\line(1,0){8}}
\put(20,150){\circle*{8}}
\put(17,156){\scriptsize{-2}}
\put(24,150){\line(1,0){8}}
\put(36,150){\circle*{8}}
\put(33,156){\scriptsize{-2}}
\put(40,150){\line(1,0){8}}
\put(56,150){\circle{16}}
\put(54,147){\footnotesize{1}}
\put(64,150){\line(1,0){8}}
\put(72,148){\fbox{$C_1$}}
\put(74,133){\scriptsize{-14}}
\put(82,159){\line(0,1){8}}
\put(82,175){\circle{16}}
\put(80,172){\footnotesize{2}}
\put(82,183){\line(0,1){8}}
\put(82,195){\circle*{8}}
\put(88,193){\scriptsize{-2}}
\put(82,199){\line(0,1){8}}
\put(82,211){\circle*{8}}
\put(88,209){\scriptsize{-2}}
\put(82,215){\line(0,1){8}}
\put(82,227){\circle*{8}}
\put(88,225){\scriptsize{-3}}
\put(82,231){\line(0,1){8}}
\put(64,150){\line(1,0){8}}
\put(82,247){\circle{16}}
\put(80,244){\footnotesize{3}}
\put(82,255){\line(0,1){8}}
\put(82,267){\circle*{8}}
\put(88,264){\scriptsize{-4}}
\put(82,271){\line(0,1){8}}
\put(82,283){\circle*{8}}
\put(88,280){\scriptsize{-2}}
\put(92,150){\line(1,0){8}}
\put(108,150){\circle{16}}
\put(106,147){\footnotesize{4}}
\put(116,150){\line(1,0){8}}
\put(128,150){\circle*{8}}
\put(125,156){\scriptsize{-2}}
\put(132,150){\line(1,0){8}}
\put(144,150){\circle*{8}}
\put(141,156){\scriptsize{-2}}
\put(148,150){\line(1,0){8}}
\put(160,150){\circle*{8}}
\put(157,156){\scriptsize{-3}}
\put(164,150){\line(1,0){8}}
\put(180,150){\circle{16}}
\put(178,147){\footnotesize{5}}
\put(188,150){\line(1,0){8}}
\put(200,150){\circle*{8}}
\put(197,156){\scriptsize{-4}}
\put(204,150){\line(1,0){8}}
\put(216,150){\circle*{8}}
\put(213,156){\scriptsize{-2}}
\put(220,150){\line(1,0){8}}
\put(236,150){\circle{16}}
\put(234,147){\footnotesize{6}}
\put(244,150){\line(1,0){8}}
\put(256,150){\circle*{8}}
\put(253,156){\scriptsize{-7}}
\put(260,150){\line(1,0){8}}
\put(276,150){\circle{16}}
\put(274,147){\footnotesize{7}}
\put(284,150){\line(1,0){8}}
\put(296,150){\circle*{8}}
\put(293,156){\scriptsize{-2}}
\put(300,150){\line(1,0){8}}
\put(312,150){\circle*{8}}
\put(309,156){\scriptsize{-4}}
\put(316,150){\line(1,0){8}}
\put(332,150){\circle{16}}
\put(330,147){\footnotesize{8}}
\put(332,142){\line(0,-1){8}}
\put(332,130){\circle*{8}}
\put(338,128){\scriptsize{-3}}
\put(332,126){\line(0,-1){8}}
\put(332,114){\circle*{8}}
\put(338,112){\scriptsize{-2}}
\put(332,110){\line(0,-1){8}}
\put(332,98){\circle*{8}}
\put(338,96){\scriptsize{-2}}
\put(332,94){\line(0,-1){8}}
\put(332,78){\circle{16}}
\put(330,75){\footnotesize{9}}
\put(332,70){\line(0,-1){8}}
\put(323,51){\fbox{$C_2$}}
\put(346,50){\scriptsize{-14}}
\put(332,46){\line(0,-1){8}}
\put(332,30){\circle{16}}
\put(327,27){\footnotesize{10}}
\put(332,22){\line(0,-1){8}}
\put(332,10){\circle*{8}}
\put(338,8){\scriptsize{-2}}
\put(332,6){\line(0,-1){8}}
\put(332,-6){\circle*{8}}
\put(338,-8){\scriptsize{-2}}
\put(328,-6){\line(-1,0){14}}
\put(310,-6){\circle*{8}}
\put(307,-18){\scriptsize{-3}}
\put(306,-6){\line(-1,0){8}}
\put(290,-6){\circle{16}}
\put(286,-9){\footnotesize{11}}
\put(282,-6){\line(-1,0){8}}
\put(270,-6){\circle*{8}}
\put(267,-18){\scriptsize{-4}}
\put(266,-6){\line(-1,0){8}}
\put(254,-6){\circle*{8}}
\put(251,-18){\scriptsize{-2}}
\put(250,-6){\line(-1,0){8}}
\put(234,-6){\circle{16}}
\put(229,-9){\footnotesize{12}}
\put(226,-6){\line(-1,0){8}}
\put(214,-6){\circle*{8}}
\put(211,-18){\scriptsize{-7}}
\put(214,-2){\line(0,1){8}}
\put(214,14){\circle{16}}
\put(209,11){\footnotesize{13}}
\put(214,22){\line(0,1){8}}
\put(214,34){\circle*{8}}
\put(201,32){\scriptsize{-2}}
\put(214,38){\line(0,1){12}}
\put(214,54){\circle*{8}}
\put(201,52){\scriptsize{-4}}
\put(218,54){\line(1,0){8}}
\put(234,54){\circle{16}}
\put(229,51){\footnotesize{14}}
\put(242,54){\line(1,0){8}}
\put(254,54){\circle*{8}}
\put(251,60){\scriptsize{-3}}
\put(258,54){\line(1,0){8}}
\put(270,54){\circle*{8}}
\put(267,60){\scriptsize{-2}}
\put(274,54){\line(1,0){8}}
\put(286,54){\circle*{8}}
\put(283,60){\scriptsize{-2}}
\put(290,54){\line(1,0){8}}
\put(306,54){\circle{16}}
\put(301,51){\footnotesize{15}}
\put(314,54){\line(1,0){9}}
\end{picture}
\end{equation*}
\end{center}

{\it The notations are as in proposition \ref{extr1}.
Every curve with discrepancy $0$ has a number between 1 and 15.
Their self-intersection indexes (on a minimal resolution)
are equal to $-1$.}

\begin{corollary}\label{extr3}
Let us consider the case $\mathrm{(I_2^2)}$ of theorem \ref{ldmain}.
Then the extraction of all exceptional curves with discrepancies
$0$ for the pair $(\overline{S},\frac67C_1+\frac67C_2)$ is shown
in the following figure.
\end{corollary}

\begin{center}
\begin{equation*}
\begin{picture}(305,300)(0,0)
\put(4,160){\circle*{8}}
\put(1,166){\scriptsize{-3}}
\put(8,160){\line(1,0){8}}
\put(20,160){\circle*{8}}
\put(17,166){\scriptsize{-2}}
\put(24,160){\line(1,0){8}}
\put(36,160){\circle*{8}}
\put(33,166){\scriptsize{-2}}
\put(40,160){\line(1,0){8}}
\put(56,160){\circle{16}}
\put(54,157){\footnotesize{1}}
\put(64,160){\line(1,0){8}}
\put(72,158){\fbox{$C_1$}}
\put(96,158){\scriptsize{-14}}
\put(82,169){\line(0,1){8}}
\put(82,185){\circle{16}}
\put(80,182){\footnotesize{4}}
\put(82,193){\line(0,1){8}}
\put(82,205){\circle*{8}}
\put(88,203){\scriptsize{-2}}
\put(82,209){\line(0,1){8}}
\put(82,221){\circle*{8}}
\put(88,219){\scriptsize{-2}}
\put(82,225){\line(0,1){8}}
\put(82,237){\circle*{8}}
\put(88,235){\scriptsize{-3}}
\put(82,241){\line(0,1){8}}
\put(82,257){\circle{16}}
\put(80,254){\footnotesize{5}}
\put(82,265){\line(0,1){8}}
\put(82,277){\circle*{8}}
\put(88,274){\scriptsize{-4}}
\put(82,281){\line(0,1){8}}
\put(82,293){\circle*{8}}
\put(79,299){\scriptsize{-2}}
\put(86,293){\line(1,0){8}}
\put(102,293){\circle{16}}
\put(100,290){\footnotesize{6}}
\put(110,293){\line(1,0){8}}
\put(122,293){\circle*{8}}
\put(119,299){\scriptsize{-7}}
\put(126,293){\line(1,0){8}}
\put(142,293){\circle{16}}
\put(140,290){\footnotesize{7}}
\put(150,293){\line(1,0){8}}
\put(162,293){\circle*{8}}
\put(159,299){\scriptsize{-2}}
\put(162,169){\line(0,1){8}}
\put(162,185){\circle{16}}
\put(160,182){\footnotesize{9}}
\put(162,193){\line(0,1){8}}
\put(162,205){\circle*{8}}
\put(168,203){\scriptsize{-2}}
\put(162,209){\line(0,1){8}}
\put(162,221){\circle*{8}}
\put(168,219){\scriptsize{-2}}
\put(162,225){\line(0,1){8}}
\put(162,237){\circle*{8}}
\put(168,235){\scriptsize{-3}}
\put(162,241){\line(0,1){8}}
\put(162,257){\circle{16}}
\put(160,254){\footnotesize{8}}
\put(162,265){\line(0,1){8}}
\put(162,277){\circle*{8}}
\put(168,274){\scriptsize{-4}}
\put(162,281){\line(0,1){8}}
\put(152,158){\fbox{$C_2$}}
\put(136,158){\scriptsize{-14}}
\put(172,160){\line(1,0){8}}
\put(188,160){\circle{16}}
\put(186,157){\footnotesize{2}}
\put(196,160){\line(1,0){8}}
\put(208,160){\circle*{8}}
\put(205,166){\scriptsize{-2}}
\put(212,160){\line(1,0){8}}
\put(224,160){\circle*{8}}
\put(223,166){\scriptsize{-2}}
\put(228,160){\line(1,0){8}}
\put(240,160){\circle*{8}}
\put(237,166){\scriptsize{-3}}
\put(244,160){\line(1,0){8}}
\put(260,160){\circle{16}}
\put(258,157){\footnotesize{3}}
\put(268,160){\line(1,0){8}}
\put(280,160){\circle*{8}}
\put(277,166){\scriptsize{-4}}
\put(284,160){\line(1,0){8}}
\put(296,160){\circle*{8}}
\put(293,166){\scriptsize{-2}}
\put(82,153){\line(0,-1){8}}
\put(82,137){\circle{16}}
\put(78,134){\footnotesize{10}}
\put(82,129){\line(0,-1){8}}
\put(82,117){\circle*{8}}
\put(88,115){\scriptsize{-2}}
\put(82,113){\line(0,-1){8}}
\put(82,101){\circle*{8}}
\put(88,99){\scriptsize{-2}}
\put(82,97){\line(0,-1){8}}
\put(82,85){\circle*{8}}
\put(88,83){\scriptsize{-3}}
\put(82,81){\line(0,-1){8}}
\put(82,65){\circle{16}}
\put(78,63){\footnotesize{11}}
\put(82,57){\line(0,-1){8}}
\put(82,45){\circle*{8}}
\put(88,43){\scriptsize{-4}}
\put(82,41){\line(0,-1){8}}
\put(82,29){\circle*{8}}
\put(79,17){\scriptsize{-2}}
\put(86,29){\line(1,0){8}}
\put(102,29){\circle{16}}
\put(97,26){\footnotesize{12}}
\put(110,29){\line(1,0){8}}
\put(122,29){\circle*{8}}
\put(119,17){\scriptsize{-7}}
\put(126,29){\line(1,0){8}}
\put(142,29){\circle{16}}
\put(137,26){\footnotesize{13}}
\put(150,29){\line(1,0){8}}
\put(162,29){\circle*{8}}
\put(159,17){\scriptsize{-2}}
\put(162,153){\line(0,-1){8}}
\put(162,137){\circle{16}}
\put(157,134){\footnotesize{15}}
\put(162,129){\line(0,-1){8}}
\put(162,117){\circle*{8}}
\put(168,115){\scriptsize{-2}}
\put(162,113){\line(0,-1){8}}
\put(162,101){\circle*{8}}
\put(168,99){\scriptsize{-2}}
\put(162,97){\line(0,-1){8}}
\put(162,85){\circle*{8}}
\put(168,83){\scriptsize{-3}}
\put(162,81){\line(0,-1){8}}
\put(162,65){\circle{16}}
\put(157,62){\footnotesize{14}}
\put(162,57){\line(0,-1){8}}
\put(162,45){\circle*{8}}
\put(168,43){\scriptsize{-4}}
\put(162,41){\line(0,-1){8}}
\end{picture}
\end{equation*}
\end{center}

{\it The notations are as in proposition \ref{extr1}.
Every curve with discrepancy $0$ has a number between 1 and 15.
Their self-intersection indexes (on a minimal resolution)
are equal to $-1$.}

\subsection*{Classification of log Enriques surfaces in the case $\mathbf{(A_2^6)}$}

Let $S$ be a log Enriques surface with $\delta=2$.
Assume that its model has the type
$\mathrm{(A_2^6)}$ (see theorem \ref{ldmain}). Then
$S$ can be constructed by the following way:
at first we extract some set $\TT$ of exceptional curves with discrepancy 0. After
it we contract the proper transforms
of $C_1$ and $C_2$. So, the classification of log Enriques surfaces with $\delta=2$
is reduced to the description of sets $\TT$.

\begin{definition}\label{def}
Let us extract some set $\TT$ of exceptional curves with discrepancy 0. Then
$\TT=\TT_1\cup\TT_2\cup\TT_3$, where
$\TT_1\subset\{1,2,3\}$, $\TT_2\subset\{4,5,6,7,8,9\}$ and
$\TT_3\subset\{10,11,12,13,14,15\}$.
Every number denotes the corresponding exceptional curve with discrepancy
0 (see the figure in corollary \ref{extr2}).
Let $T_i=|\TT_i|$ be a number of elements of set $\TT_i$ and
$\overline{T_i}=\max\{t\mid t\in\TT_i\}$, $\underline{T_i}=\min\{t\mid t\in\TT_i\}$.
\end{definition}

In the following first classification theorem the set
$\TT_3$ is considered up to symmetry.
For example, if $10\notin \TT_3$ then $15\notin \TT_3$
\begin{theorem}\label{main1}
In the case $\mathrm{(A_2^6)}$ the set $\TT$ must satisfy the condition $(1)$ and
be one of the following sets:
\begin{enumerate}
\item Always $T_1+T_2\ge 1$ and $T_3\ge 1$. If $T_1=0$ then $\TT_2\ne \{9\}$.
\item Let $T_2=0$ and $T_3=1$. Then $10\in \TT_3$.
\item Let $T_2\ge 1$ and $T_3=1$. If $\overline T_2\le 6$ then $10\in \TT_3$.
If $\overline T_2=7$ then $\{10,11\}\cap \TT_3\ne \emptyset$.
If $\overline T_2=8\ \text{or}\ 9$ then $\TT_3\ne \emptyset$ is arbitrary.
\item Let $T_2=0$ and $T_3=2$. Then either $10\in \TT_3$, or
$\TT_3=\{11,14\}$ and $\TT_1\ne \{3\}$.
\item Let $T_2\ge 1$ and $T_3=2$. Then either $10\in \TT_3$, or
$\TT_3=\{11,14\}$, or $\TT_3=\{11,13\}$ and $\overline{T}_2\ge 5$, or
$\TT_3=\{11,12\}$ and $\{7,8,9\}\cap \TT_2\ne \emptyset$, or
$\TT_3=\{12,13\}$ and $\{8,9\}\cap \TT_2\ne \emptyset$.
\item Let $T_2=0$ and $T_3=3$. Then either $10\in \TT_3$, or
$\TT_3=\{11,12,14\}$ and $\TT_1\ne \{3\}$.
\item Let $T_2\ge 1$ and $T_3=3$. Then $10\in \TT_3$, or
$\TT_3=\{11,12,14\}$, or $\TT_3=\{11,12,13\}$ and $\overline{T}_2\ge 5$.
\item Let $T_2=0$ and $T_4=4$. Then either $10\in \TT_3$, or
$\TT_3=\{11,12,13,14\}$ and $\TT_1\ne \{3\}$.
\item If either $T_2\ge 1$ and $T_3\ge 4$, or $T_3\ge 5$ then
$\TT$ is arbitrary $($of course, taking into account the condition $(1))$.
\end{enumerate}
\begin{proof} Let $f\colon Y\to X$ be an extraction of some set $\TT$.
The proper transforms of $C_i$ are denoted by $\widetilde{C}_i$. We have
$C_1^2=\frac16$ and $C_2^2=6$.
The obvious requirement for the set $\TT$ is
$\widetilde{C}_i^2<0$. Therefore $T_1+T_2\ge 1$ and $T_3\ge 1$. Hence, always
$\widetilde{C}_1^2<0$, except one case
$T_1=0$ and $\TT_2=\{9\}$. For this case $\widetilde{C}_1^2=0$.
\par
The remaining cases are written taking into account the requirements
$\widetilde{C}_2^2<0$. We can contract $\widetilde{C}_1$ and
$\widetilde{C}_2$ except one case, which appears in the conditions 4, 6, 8.
\par
For this case
$T_2=0$, $2\le T_3\le 4$, $\TT_1=\{3\}$, $10\notin \TT_3 \supset \{11,14\}$.
Consider the minimal resolution. Then our configuration is illustrated in
the following figure.\\
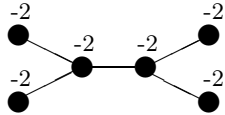
\begin{center}
\begin{picture}(80,32)(0,0)
\put(5,30){\circle*{8}}
\put(2,36){\scriptsize{-2}}
\put(5,30){\line(2,-1){24}}
\put(5,5){\circle*{8}}
\put(2,11){\scriptsize{-2}}
\put(5,6){\line(2,1){24}}
\put(29,18){\circle*{8}}
\put(26,24){\scriptsize{-2}}
\put(33,18){\line(1,0){16}}
\put(53,18){\circle*{8}}
\put(50,24){\scriptsize{-2}}
\put(53,18){\line(2,1){24}}
\put(50,18){\line(2,-1){24}}
\put(77,30){\circle*{8}}
\put(75,36){\scriptsize{-2}}
\put(77,5){\circle*{8}}
\put(75,11){\scriptsize{-2}}
\end{picture}
\end{center}
Indeed, the determinant of intersection matrix is equal to 0, although
$\widetilde{C}_1^2=\widetilde{C}_2^2=-1$.
\end{proof}
\end{theorem}

\subsection*{Classification of log Enriques surfaces in the case $\mathbf{{(I_2^2)}}$}

We have the same notations as in definition
\ref{def}. The corresponding figure is given in corollary \ref{extr3}.
The case $T_2=0$, $T_3\ge 1$ is symmetric to the case
$T_2\ge 1$, $T_3=0$ and therefore it isn't considered.
\begin{theorem}\label{main2}
In the case $\mathrm{(I_2^2)}$ the set $\TT$ must satisfy the condition $(1)$ and
be one of the following sets:
\begin{enumerate}
\item $T_2+T_3\ge 1$.
\item Let $T_3=0$ and $T_1=0$. Then $\underline{T_2}\le 5$ and $9\in \TT_2$.
\item Let $T_3=0$ and $T_1=1$. Then either $\TT_1=\{1\}$ and $9\in \TT_2$, or
$\TT_1=\{2\}$ and $\underline{T_2}\le 5$, or $\TT_1=\{3\}$ and
$4\in \TT_2 \cap \{8,9\}\ne \emptyset$, or $\TT_1=\{3\}$ and
$ \{5,9\}\subset \TT_2$.
\item Let $T_3=0$ and $T_1=2$. Then either $\TT_1=\{1,2\}$, or $\TT_1=\{1,3\}$ and
$\TT_2 \cap \{8,9\}$, or $\TT_1=\{2,3\}$ and $\underline{T_2}\le 5$.
\item Let $T_3=0$ and $T_1=3$. Then $\TT$ is arbitrary.
\item Let $T_2,T_3\ge 1$ and $T_1=0$. Then the set $\TT$ must satisfy the following
two conditions $\Upsilon_1$ and $\Upsilon_2$.\\
\underline{Condition $\Upsilon_1$}.
$(\underline{T_2}\le 5, \TT_3\ \text{is arbitrary})$,
or $(\underline{T_2}=6, \underline{T_3}\le 14)$,
or $(7\le \underline{T_2}\le 8, \underline{T_3}\le 12)$,
or $(\underline{T_2}=9, \underline{T_3}\le 11)$.\\
\underline{Condition $\Upsilon_2$}. $(\overline{T_2}=4, T_3=15)$,
or $(5\le\overline{T_2}\le 7, \overline{T_3}\ge 14)$,
or $(\overline{T_2}=8, \overline{T_3}\ge 11)$,
or $(\overline{T_2}=9, \TT_3\  \text{is arbitrary})$.
\item Let $T_2,T_3\ge 1$ and $\TT_1=\{1\}$. Then the set $\TT$ must satisfy
the condition $\Upsilon_2$.
\item Let $T_2,T_3\ge 1$ and $2\in \TT_1$. Then the set $\TT$ must satisfy
the condition $\Upsilon_1$.
\item Let $T_2,T_3\ge 1$ and $\TT_1=\{3\}$. Then the set $\TT$ must satisfy
the condition $\Upsilon_1$ and the next condition $\Upsilon_3$.\\
\underline{Condition $\Upsilon_3$}. $(\overline{T_2}=4, \overline{T_3}\ge 14)$,
or $(5\le\overline{T_2}\le 6, \overline{T_3}\ge 13)$,
or $(\overline{T_2}=7, \overline{T_3}\ge 11)$,
or $(\overline{T_2}\ge 8, \TT_3 \ \text{is arbitrary})$.
\item Let $T_2,T_3\ge 1$ and $\{1,2\}\subset\TT_1$. Then $\TT$ is arbitrary.
\item Let $T_2,T_3\ge 1$ and $\TT_1=\{1,3\}$. Then the set $\TT$ must satisfy
the condition $\Upsilon_3$.\\
\end{enumerate}
\begin{proof} Let $f\colon Y\to X$ be an extraction of some set $\TT$.
The proper transforms of $C_i$ are denoted by $\widetilde{C}_i$. We have
$C_1^2=\frac32$ and $C_2^2=\frac83$.
If $T_2=T_3=0$ then $\widetilde{C}_1^2<0$ and $\widetilde{C}_2^2<0$ in the case
$\{1,2\}\subset\TT_1$ only.
Consider the minimal resolution. Then our configuration is illustrated in
the following figure.\\
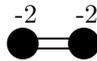
\begin{center}
\begin{picture}(36,10)(0,0)
\put(7,6){\circle*{12}}
\put(4,14){\footnotesize{-2}}
\put(7,4){\line(1,0){20}}
\put(30,6){\circle*{12}}
\put(27,14){\footnotesize{-2}}
\put(7,8){\line(1,0){20}}
\end{picture}
\end{center}
The determinant of intersection matrix is equal to 0 and therefore
always $T_2+T_3\ge 1$.
The rest cases are written taking into account the requirements
$\widetilde{C}_1^2<0$ and
$\widetilde{C}_2^2<0$.
We can contract
$\widetilde{C}_1$ and
$\widetilde{C}_2$ except one case:
$T_3=0$, $\TT_1=\{3\}$, $\underline{T_2}=5$,
$\overline{T_2}=8$.
Consider the minimal resolution. Then our configuration is illustrated in
the figure of theorem \ref{main1} proof.
It appears in the case (3).
\end{proof}
\end{theorem}

By theorems \ref{main1} and \ref{main2} we get the following corollary.
\begin{corollary} Let $S$ be a log Enriques surface with $\delta=2$. Then
$1\le \rho(S)\le 14$.
\end{corollary}

\subsection*{Log Enriques surfaces and $K3$ surfaces}

Let $S$ be a log Enriques surface with $\delta=2$.
Consider its canonical cover
$\varphi\colon \widehat{S}=\Spec_{\OO_S}(\oplus^{I-1}_{i=0}\OO_S(-iK_S))\to S$, where
$I$ is an index of $S$. Since $I=7$ (see corollary \ref{index}) then
$\widehat{S}$ is a $K3$ surface with at worst Du Val singularities
\cite{Bl}, \cite{Z} and
\begin{enumerate}
\item $\varphi$ is cyclic Galous cover of degree 7, which is
etale over $S\backslash\Sing S$.
\item There exists a generator $g$ of $\Gal(\widehat{S}/S)\cong \ZZ_7$
such that $g^*\omega_{\widehat{S}}=\varepsilon_7\omega_{\widehat{S}}$, where
$\varepsilon_7=\exp(2\pi \sqrt{-1}/7)$ is a primitive root
and $\omega_{\widehat{S}}$ is a
nowhere vanishing regular
2-form on $\widehat{S}$.
\end{enumerate}
Let $\Delta(S)$ be an exceptional set of $\chi$, where $\chi$ is a minimal resolution
of
$\widehat{S}$. Then $\Delta(S)$ is a disconnected sum of divisors of Dynkin's type
$\AAA_i$, $\DDD_j$, $\EEE_k$.
So $\Delta(S)=(\oplus \AAA_{\alpha})\oplus(\oplus \DDD_{\beta})
\oplus(\oplus \EEE_{\gamma})$. Let us define
$\rank\Delta(S)=\sum\alpha+\sum\beta+\sum\gamma$.

\begin{theorem} Let $S$ be a log Enriques surface with $\delta=2$. Then
$\rank\Delta(S)+\rho(S)=16$.
\begin{proof} Let $\TT=\{1,2,\ldots,15\}$.
The surface corresponding to this set is denoted by $S'$. In the cases
$\mathrm{(A^6_2)}$ and
$\mathrm{(I^2_2)}$ the set $\Delta(S')$ has the same type
$\AAA_1\oplus\AAA_1$, i.e.
$\rank\Delta(S')+\rho(S')=2+14=16$. Every surface
$S$ can be obtained from $S'$ by the contractions of some curves in
$\TT$. If we contract any curve in $\TT$ then
$\rank \Delta$ is increased by 1.
Therefore the required equality is reserved.
\end{proof}
\end{theorem}

\begin{example}
Let $\TT=\{8,12\}$ then $\rho(S)=1$ and $\Delta(S)=\EEE_7\oplus\EEE_8$ in both cases
$\mathrm{(A^6_2)}$ and $\mathrm{(I^2_2)}$.
\end{example}

\begin{remark}
Other approach to the classification was given in \cite{Z2}.
In particular, see \cite[\S 5]{Z2} in the case of surfaces with index 7.
\end{remark}

\end{document}